\documentclass[]{article}
\usepackage{latexsym}
\usepackage{amssymb,latexsym,amsmath,enumerate,verbatim,amsfonts}
\newtheorem{definition}{Definition}[section]
\newtheorem{theorem}{Theorem}[section]
\newtheorem{corollary}[theorem]{Corollary} 

\newtheorem{lemma}{Lemma}[section]

\newtheorem{remark}{Remark}[section]
\numberwithin{equation}{section}
\title{Positive definiteness of a class of cyclic symmetric tensors}
\author{Yisheng Song
\thanks{The author’s work was supported by the National Natural Science Foundation of P.R. China (Grant No.12171064), by The team project of innovation leading talent in chongqing (No.CQYC20210309536) and by the Foundation of Chongqing Normal university (20XLB009).}\\
{\small School of Mathematical Sciences, Chongqing Normal University,}\\{\small Chongqing 401331 P.R. China.
	Email: yisheng.song@cqnu.edu.cn} }
\date{ }
\begin{document}

\maketitle

\begin{abstract}
For a 4th order 3-dimensional cyclic symmetric tensor, a  sufficient and necessary  condition is bulit for  its positive semi-definiteness. A sufficient and necessary  condition of positive definiteness is showed for a 4th order $n$-dimensional  symmetric tensor. With the help of such a condition,   the positive definiteness of a class of 4th order 3-dimensional cyclic symmetric  tensors is given. Moreover, the positive definiteness of a class of non-cyclic symmetric  tensors is showed also.  By applying these conclusions, several (strict) inequalities are erected for ternary quartic homogeneous polynomials.
\end{abstract}

\section{Introduction}
As we all know,  a fundamental and challenging problem in dealing with tensors is to determine its positive definiteness, which is  an NP-hard problem in general even when the order is $4$. The notion of positive definiteness for a symmetric tensor was first used   by Qi \cite{Q2005}.  The  vacuum stability of the Higgs scalar potential  model in High-Energy Physics is one of the most direct applications of positive definiteness of 4th order tensors \cite{K2016,SQ2024,S2023,S20231}.  In order to verifying the positive definiteness of higher order tensors (higher homogeneous polynomial), the notion of cyclic symmetric group was introduced, which may trace  back to ones of  Refs. Stanley \cite{S1979,S2007}, Yuan \cite{Y2007}, Gao-Geroldinger \cite{GG2006}, Harris \cite{HW2006}. However, even if a tensor (or homogeneous polynomial) is cyclic symmetric with the order $4$ and dimension $3$, a prefect criterion of positive definiteness was not obtained untill now.

Since $2-$dimensional symmetric tensor can be turn into a univariate polynomial and the positivity condition of a quartic univariate polynomial has had prefect identifying criterion, so  the positive definiteness of a 4th order $2-$dimensional symmetric tensor can be well determine and has analytic criterion. Such a criterion dates back to 1922 (Rees \cite{R1922}).  Subsequently,  many academics devoted themself to make sure its identifying conditions  such as ones of Refs. Gadem-Li \cite{GL1964},  Ku \cite{K1965},  Jury-Mansour \cite{JM1981} and  Lazard \cite{L1988}.  Untill to 2005, Wang-Qi \cite{WQ2005}  perfectly gave  analytic necessary and sufficient conditions for checking the positive definiteness of a 4th order $2-$dimensional symmetric tensor. In 1996,  Hasan-Hasa \cite{HH1996} claimed that they provided a necessary and sufficient condition of the positive definiteness of higher order tensors (higher then 4th onder)  without the discriminant. However, in 1998, Fu \cite{F1998} pointed out that Hasan-Hasan's results are sufficient only.  Recently, Guo\cite{G2021} showed a new necessary and sufficient condition without the discriminant for  a 4th order $2-$dimensional symmetric tensor. Very recently, Qi-Song-Zhang \cite{QSZ2022} gave a analytic necessary and sufficient condition, which treats  cubic and quartic coefficients in a symmetric way, and  explicitly states  a necessary condition.

In 2005, Qi \cite{Q2005} presented to verify the positive definitiveness of  a even order symmetric tensor by means of  the sign of its all H-(Z-)eigenvalue.  Subsequently, many scholars d committed  themself to solve or compute H-(Z-)eigenvalue of a tensor.  Ni-Qi-Wang \cite{NQW2008}  provided  a method of  checking positive definiteness of a 4th order 3 dimensional tensor by computing its smallest eigenvalue. Song \cite{S2021} and Song-Qi \cite{SQ2021} proved several  sufficient conditions and necessary conditionsof its positive definiteness of a 4th order  symmetric tensor.   Ng-Qi-Zhou \cite{NQZ2009} gave   an algorithm of the largest eigenvalue of a nonnegative tensor.  Cui-Dai-Nie \cite{CDN2014} presented a way of computing all real eigenvalues of a symmetric tensor.  Zhang-Qi \cite{ZQ2012} gave linear convergence of  the largest eigenvalue of tensors. Hu-Huang-Qi \cite{HHQ2013,HHQ2014} presented some properities of H-(Z-)eigenvalue of tensors and their calculation methods. For more details about  H-(Z-)eigenvalues, see Refs. \cite{HHLQ2013, HLQS2013,LWZ2014,CCW2016,QL2018, QCC2018, DQW2013,KM2011,ZQZ2014,YY2010,YY2011,KM2014,CHZ2016,KSB2015,LZI2010,NQ2015,NZ2018,ZCQ2013,ZWQA2018,LN2015,SQ2016,CW2017,ZQX2012,MDW2023,LQY2015,HJ2016,ZQW2013}  and other references not cited here.

In this paper,  we first dicuss  necessary and sufficient conditions of positive (semi-)definiteness of a class of  symmetric tensors (Theorems \ref{thm:41} and \ref{thm:42}).  For a cyclic symmetric tensor of orrder 4 and dimension 3, the conditions determining its positive definiteness  are given (Theorems \ref{thm:43} and \ref{thm:44}).  The condition checking the positive definiteness of a class of non-cyclic symmetric tensors is presented also(Theorem \ref{thm:46}).   Furthermore,   several  (strict) inequalities of  quartic ternary homogeneous polynomial  are built (Corollaries \ref{cor:47} and \ref{cor:48}) with the help of the above conclusions about tensors.

\section{Positive definiteness of 4th order 2-dimensional symmetric tensors}

An $d$th order $n$ dimensional   tensor $\mathcal{T}=(t_{i_1i_2\cdots i_d})$ is said to be {\bf symmetric}  if $$t_{i_1i_2\cdots i_d}=t_{j_1j_2\cdots j_d}$$ whenever $(i_1, i_2,\cdots,i_d)$ is a permutation of $(j_1, j_2,\cdots,j_d)$.

\begin{definition} \em Let $\mathcal{T}=(t_{i_1i_2\cdots i_d})$  be an
	$d$th order $n$ dimensional symmetric tensor with even number $d$.  $\mathcal{T}$ is called 
	\begin{itemize}
		\item[(i)] {\bf positive semi-definite} (\cite{Q2005}) if   in the Euclidean space $ \mathbb{R}^n$, its associated Homogeneous polynomial $$\mathcal{T}x^d=\sum\limits_{i_1,i_2,\cdots,i_m=1}^nt_{i_1i_2\cdots
			i_m}x_{i_1}x_{i_2}\cdots
		x_{i_m}\geq0;$$ 
		\item[(ii)] {\bf positive definite} (\cite{Q2005}) if   $\mathcal{T}x^d>0$ for all $x\in \mathbb{R}^n\setminus\{0\}$.
	\end{itemize} 
\end{definition}

Let $\mathcal{T}=(t_{ijkl})$ be a 4th-order 2-dimensional symmetric tensor. Then for $x=(x_1,x_2)^\top,$
\begin{equation}\label{eq:f}
	Tx^4=t_{1111} x_1^4+4t_{1112} x_1^3x_2+6t_{1122} x_1^2x_2^2+4t_{1222}x_1x_2^3+t_{2222}x_2^4.
\end{equation}
Let $$\begin{aligned}
	\Delta	=&4\times 12^3(t_{1111}t_{2222}-4t_{1112}t_{1222}+3t_{1122}^2)^3\\
	&-72^2\times 6^2(t_{1111}t_{1122}t_{2222}+ 2t_{1112}t_{1122}t_{1222}- t_{1122}^3-t_{1111}t_{1222}^2- t_{1112}^2t_{2222})^2\\
	=& 4\times 12^3(\eta^3-27\chi^2),
\end{aligned}$$
where $$\begin{aligned}\eta=&t_{1111}t_{2222}-4t_{1112}t_{1222}+3t_{1122}^2,\\
	\chi=&t_{1111}t_{1122}t_{2222}+2t_{1112}t_{1122}t_{1222}-t_{1122}^3-t_{1111}t_{1222}^2-t_{1112}^2t_{2222}.
\end{aligned}$$
and hence, the sign of $\Delta$ is the same as one of  $(\eta^3-27\chi^2)$. 
\begin{lemma}[\cite{SQ2024,QSZ2022}]\label{lem:21} A 4th-order 2-dimensional symmetric tensor $\mathcal{T}=(t_{ijkl})$ is positive definite if and only if
	$$\begin{cases}
	\eta^3-27\chi^2=0,\ \ t_{1112}\sqrt{t_{2222}}=t_{1222}\sqrt{t_{1111}},\\
		2t_{1112}^2+t_{1111}\sqrt{t_{1111}t_{2222}}=3t_{1111}t_{1122}<3t_{1111}\sqrt{t_{1111}t_{2222}};\\
		\eta^3-27\chi^2>0,\\
		|t_{1112}\sqrt{t_{2222}}-t_{1222}\sqrt{t_{1111}}|\leq \sqrt{6t_{1111}t_{1122}t_{2222}+2\sqrt{(t_{1111}t_{2222})^3}},\\
		(i) \ -\sqrt{t_{1111}t_{2222}}< 3t_{1122}\leq 3\sqrt{t_{1111}t_{2222}};\\
		(ii)\  t_{1122} >\sqrt{t_{1111}t_{2222}}\ \mbox{ and } \\ |t_{1112}\sqrt{t_{2222}}+t_{1222}\sqrt{t_{1111}}|\leq \sqrt{6t_{1111}t_{1122}t_{2222}-2\sqrt{(t_{1111}t_{2222})^3}}.
	\end{cases}$$
	A 4th-order 2-dimensional symmetric tensor $\mathcal{T}=(t_{ijkl})$ is positive semidefinite if and only if
	$$\begin{cases}
		\eta^3-27\chi^2\ge0,\\
		|t_{1112}\sqrt{t_{2222}}-t_{1222}\sqrt{t_{1111}}|\leq \sqrt{6t_{1111}t_{1122}t_{2222}+2\sqrt{(t_{1111}t_{2222})^3}},\\
		(i) \ -\sqrt{t_{1111}t_{2222}}\leq 3t_{1122}\leq 3\sqrt{t_{1111}t_{2222}};\\
		(ii) 	\ t_{1122} >\sqrt{t_{1111}t_{2222}}\mbox{ and }\\
		|t_{1112}\sqrt{t_{2222}}+t_{1222}\sqrt{t_{1111}}|\leq \sqrt{6t_{1111}t_{1122}t_{2222}-2\sqrt{(t_{1111}t_{2222})^3}}.
	\end{cases}$$
\end{lemma}
\begin{corollary}\label{cor:21} Let $\mathcal{T}=(t_{ijkl})$  be a 4th-order 2-dimensional symmetric tensor with its entires $|t_{ijkl}|\leq1$ and $t_{1111}=t_{2222}=1.$  Then \begin{itemize}
		\item[(i)]   $\mathcal{T}$ is positive semi-definite if and only if$$\begin{cases}
			-\dfrac13< t_{1122}\leq1,	(t_{1112}-t_{1222})^2\leq6t_{1122}+2,\\
			27(t_{1122}+2t_{1112}t_{1122}t_{1222}-t_{1122}^3-t_{1222}^2-t_{1112}^2)^2\le(1-4t_{1112}t_{1222}+3t_{1122}^2)^3.
		\end{cases}$$  
		\item[(ii)]   $\mathcal{T}$ is positive definite if and only if $$\begin{cases}
			\dfrac13\leq t_{1122}<1$, $2t_{1112}^2+1=3t_{1122}$, $t_{1112}=t_{1222}; \\
			-\dfrac13< t_{1122}\leq1,	(t_{1112}-t_{1222})^2\leq6t_{1122}+2,\\
			27(t_{1122}+2t_{1112}t_{1122}t_{1222}-t_{1122}^3-t_{1222}^2-t_{1112}^2)^2< (1-4t_{1112}t_{1222}+3t_{1122}^2)^3.
		\end{cases}$$  
	\end{itemize}
\end{corollary}
{\bf Proof.}
	(i) It follows from Theorem \ref{lem:21}  that $\mathcal{T}$ is positive semi-definite if and only if
	$$\eta^3-27\chi^2\ge0,\ \  |t_{1112}-t_{1222}|\leq\sqrt{6t_{1122}+2}\mbox{ and }-1\leq 3t_{1122} \leq3. $$
	
	(ii) By Theorem \ref{lem:21},  $\mathcal{T}$ is positive definite if and only if
	
	$$\aligned
	\eta^3-27\chi^2=0,& \ t_{1112}=t_{1222}\mbox{ and }1\leq2t_{1112}^2+1=3t_{1122}<3;\\
	\eta^3-27\chi^2>0,& \  |t_{1112}-t_{1222}|\leq\sqrt{6t_{1122}+2}\mbox{ and }-1< 3t_{1122} \leq3.
	\endaligned$$
	Since	$$\left. \aligned t_{1112}=t_{1222}\\  2t_{1112}^2+1=3t_{1122}\endaligned \right\} \Rightarrow \begin{cases}
		27\chi^2=27(3t_{1122}^2-t_{1122}^3-3t_{1122}-1)^2\\
		\ \ \ \ \ \ \ \ \ \ =27(t_{1122}(1+2t_{1112}^2)-t_{1122}^3-2t_{1112}^2)^2,\\
		\eta^3 =(1-4t_{1112}^2+3t_{1122}^2)^3\\
		\ \ \ \ \  =(1-2(3t_{1122}-1)+3t_{1122}^2)^3,
	\end{cases}$$ 
	we have $$27\chi^2=\eta^3=27(t_{1122}-1)^6.$$
	This  completes the proof.	\\

When $t_{1122} = 1, (t_{1112}-t_{1222})^2 \leq6t_{1122} + 2 = 8$ is obviously satisfied since $|t_{ijkl}| \leq1$ and $2t_{1112}^2+1=3t_{1122}<3$ can't hold. It is easy to prove the following conclusions

\begin{corollary}\label{cor:22}
	Let $\mathcal{T}=(t_{ijkl})$  be a 4th-order 2-dimensional symmetric tensor with its entires $|t_{ijkl}|\leq1$ and $t_{1111}=t_{2222}=t_{1122}=1.$  Then \begin{itemize}
		\item[(i)]   $\mathcal{T}$ is positive semi-definite if and only if
		$27(t_{1222}-t_{1112})^4\leq 64(1-t_{1112}t_{1222})^3;$
		\item[(ii)]   $\mathcal{T}$ is positive definite if and only if $27(t_{1222}-t_{1112})^4< 64(1-t_{1112}t_{1222})^3.$
	\end{itemize}
\end{corollary}

\begin{corollary}\label{cor:23} Let $\mathcal{T}=(t_{ijkl})$  be a 4th-order 2-dimensional symmetric tensor with its entires $|t_{ijkl}|=1$ and $t_{1111}=t_{2222}=1.$  Then \begin{itemize}
		\item[(i)]   $\mathcal{T}$ is positive semi-definite if and only if
		$t_{1122}=1;$
		\item[(ii)]   $\mathcal{T}$ is positive definite if and only if
		$t_{1122}=1\mbox{ and } t_{1112}t_{1222}=-1.$
	\end{itemize}
\end{corollary}
{\bf Proof.}
	(i) It follows from Corollaries \ref{cor:21} (i) and  \ref{cor:22} (i) that 
	$t_{1122}=1$ and either $t_{1112}t_{1222}=1$,  $$64(1-t_{1112}t_{1222})^{3}=27(t_{1112}-t_{1222})^{2}=0,$$ 
	or  $t_{1112}t_{1222}=-1,$
	$$ 64(1+1)^{3}>27(t_{1112}-t_{1222})^{2}=27\times 4.$$
	So $\mathcal{T}$ is positive semi-definite if and only if
	$t_{1122}=1.$
	
	(ii) It follows from Corollaries \ref{cor:21} (ii) and  \ref{cor:22} (ii) that 
	 $\mathcal{T}$ is positive definite if and only if $t_{1122}=1$ and
	$$t_{1112}t_{1222}=-1\mbox{ and }27(t_{1112}-t_{1222})^{2}=27\times 4<64(1+1)^{3}.$$ Then $\mathcal{T}$ is positive definite if and only if
	$t_{1122}=1\mbox{ and } t_{1112}t_{1222}=-1.$

\begin{corollary} \label{cro:24} Assume a 4th order 2-dimensional  symmetric tensor $\mathcal{T}$ is positive semi-definite. Then
	\begin{itemize}
		\item [(1)] $t_{1112}=0$ and $t_{1122}\ge0$ if $t_{1111}=0$ and $t_{2222}>0$;
		\item[(2)] $t_{1222}=0$  and $t_{1122}\ge0$ if $t_{2222}=0$ and $t_{1111}>0$;
		\item[(3)]  $t_{1112}=t_{1222}=0$  and $t_{1122}\ge0$ if $t_{2222}=t_{1111}=0$. 
	\end{itemize}
\end{corollary}
{\bf Proof.} (1) By Lemma \ref{lem:21},  we have 
	$$\aligned |t_{1112}\sqrt{t_{2222}}|\leq0\  \Longrightarrow\ t_{1112}=0;\\
	3t_{1122}\ge-\sqrt{t_{1111}t_{2222}}=0\  \Longrightarrow\ t_{1122}\ge0,\endaligned$$
	
		Similarly, the conclusion (2) is easy to be bulit.
	
	(3) Obviously, $t_{1122}\ge0$ by Lemma \ref{lem:21}.  Suppose $t_{1112}\ne0$ or $t_{1222}\ne0$.  Then  for all $x=(x_1,x_2)^\top\in\mathbb{R}^2$,
	$$\aligned \mathcal{T}x^4=&4t_{1112}x_1^3x_2+6t_{1122}x_1^2x_2^2+4t_{1222}x_1x^3_2\ge0.
	\endaligned$$
	In the meantime,  for $x_1\ne0$ and $x_2\ne0$, we have
	$$t_{1122}\ge-\dfrac23t_{1112}\dfrac{x_1}{x_2}-\dfrac23t_{1222}\dfrac{x_2}{x_1}.$$
	
	If $t_{1112}>0$, then take $x_1\to-\infty$ to yield $t_{1122}\ge+\infty;$  If $t_{1112}<0$, then take $x_1\to+\infty$ to yield $t_{1122}\ge+\infty,$ a contradiction. Therefore, $t_{1112}=0$.
	
	If $t_{1222}>0$, then take $x_2\to-\infty$ to yield $t_{1122}\ge+\infty;$  If $t_{1222}<0$, then take $x_2\to+\infty$ to yield $t_{1122}\ge+\infty,$ a contradiction. Therefore, $t_{1222}=0$.\\

Since any principal subtensor of a positive semi-definite tensor is still positive semi-definite, the following conclusion is obvious.

\begin{corollary}\label{cor:25}
	Assume a 4th order  symmetric tensor $\mathcal{T}$ is positive semi-definite. Then $$t_{iiii}=0  \  \Longrightarrow\  t_{iiij}=0\mbox{ and }t_{iijj}\ge0, \mbox{ for all }i,j.$$
\end{corollary}
\section{Inner Products and Outer Products}

For any two $d$th order $n$-dimensional tensors $\mathcal{T}=(t_{i_1i_2\cdots i_d})$ and $\mathcal{A}=(a_{i_1i_2\cdots i_d})$,  the \textbf{inner product} of $\mathcal{T}$ and $\mathcal{A}$\index{inner product},  denoted as $\langle\mathcal{T}, \mathcal{A}\rangle$, is definited as
$$\langle\mathcal{T}, \mathcal{A}\rangle=\sum_{i_2, \cdots, i_d = 1}^nt_{i_1i_2\cdots i_d}a_{i_1i_2\cdots i_d}.$$
The \textbf{Frobenius norm}\index{Frobenius norm} of $\mathcal{T}$, denoted as $\|\mathcal{T}\|_F$,  is given by
$$\|\mathcal{T}\|_F=\sqrt{\langle\mathcal{T}, \mathcal{T}\rangle}=\sqrt{\sum_{i_2, \cdots, i_d = 1}^nt_{i_1i_2\cdots i_d}^2}.$$

For any two vectors $x=(x_1,x_2,\cdots, x_n)^\top\in\mathbb{R}^n$ and $y=(y_1,y_2,\cdots, y_n)^\top\in\mathbb{R}^n$, it was already well known tha the definition of the \textbf{outer product} of $x$ and $y$\index{outer product}, (also known as \textbf{tensor product})\index{tensor product},  denoted as $x\otimes y$, is 
$$x\otimes y=(x_iy_j)=\begin{pmatrix}x_1y_1\ &\ x_1y_2\ &\ x_1y_3\ &\ \cdots\ &\ x_1y_n\\
	x_2y_1\ &\ x_2y_2\ &\ x_2y_3\ &\ \cdots\ &\ x_2y_n\\
	\vdots \ &\ \vdots\ &\ \vdots\ &\ \vdots\ &\ \vdots\\
	x_ny_1\ &\ x_ny_2\ &\ x_ny_3\ &\ \cdots\ &\ x_ny_n
\end{pmatrix}$$

Apparently,  $x\otimes x\otimes x=(x_ix_jx_k)$ and  $x\otimes x\otimes x \otimes x=(x_ix_jx_kx_l)$ are 3rd symmetric tensor and 4th symmetric tensor, respectively. Let $$x^{\otimes d}=\overbrace{x\otimes x\otimes\cdots \otimes x}^{d \text{ times}}=(x_{i_1}x_{i_2}\cdots x_{i_d}).$$ We call $x^{\otimes d}$ a \textbf{symmetric rank-one tensor}\index{symmetric rank-one tensor}. It is easy to check that
\begin{itemize}
	\item $\mathcal{T}x^d=\langle\mathcal{T}, x^{\otimes d}\rangle=\sum\limits_{i_1,i_2, \cdots, i_d = 1}^nt_{i_1i_2\cdots i_d}x_{i_1}x_{i_2}\cdots x_{i_d};$
	\item $\mathcal{T}x^{d-1}=\left(\langle\mathcal{T}_k, x^{\otimes (d-1)}\rangle\right)=\left(\aligned\sum\limits_{i_2, \cdots, i_d = 1}^nt_{1i_2\cdots i_d}x_{i_2}\cdots x_{i_d}\\
	\sum\limits_{i_2, \cdots, i_d = 1}^nt_{2i_2\cdots i_d}x_{i_2}\cdots x_{i_d}\\
	\vdots\ \ \ \ \ \ \ \ \ \ \ \ \ \ \ \ \ \ \ \ \ \ \\
	\sum\limits_{i_2, \cdots, i_d = 1}^nt_{ni_2\cdots i_d}x_{i_2}\cdots x_{i_d}
	\endaligned\right)$, \\$
	\ \mathcal{T}_k=(t_{ki_2\cdots i_d})\mbox{ for each }k=1,2,\cdots n;$
	\item $(x+y)\otimes  z=x\otimes  z+y\otimes  z$, $z\otimes (x+y)=z\otimes  x+z\otimes  y$, for all $x,y,z\in\mathbb{R}^n$;
	\item $(\alpha x)\otimes y=x\otimes(\alpha y)=\alpha x\otimes y,$ for all $\alpha\in\mathbb{R}.$
\end{itemize} Thus, we have
$$\aligned (\alpha x+\beta y)^{\otimes 3}=&((\alpha x_i+\beta y_i)(\alpha x_j+\beta y_j)(\alpha x_k+\beta y_k))\\
= &\alpha^3 x^{\otimes 3}+3\alpha^2\beta x^{\otimes 2}\otimes y+3\alpha\beta ^2x\otimes y^{\otimes 2}+\beta^3 y^{\otimes 3}
\endaligned$$ and  $$\aligned (\alpha x+\beta y)^{\otimes 4}=&((\alpha x_i+\beta y_i)(\alpha x_j+\beta y_j)(\alpha x_k+\beta y_k)(\alpha x_l+\beta y_l))\\
= &\alpha^4x^{\otimes 4}+4\alpha^3\beta x^{\otimes 3}\otimes y+6\alpha^2\beta^2 x^{\otimes 2}\otimes y^{\otimes 2}+3\alpha\beta ^3x\otimes y^{\otimes 3}+\beta^4 y^{\otimes 4}.
\endaligned$$
Therefore, for a $4$th order symmetric tensor  $\mathcal{T}=(t_{ijkl})$ and a $3$rd order symmetric tensor $\mathcal{A}=(a_{ijk})$, 
$$\aligned\mathcal{T}(\alpha x+\beta y)^4=\langle\mathcal{T}, (\alpha x+\beta y)^{\otimes 4}\rangle=&\alpha^4\mathcal{T}x^4+4\alpha^3\beta \mathcal{T}x^3 y\\
&+6\alpha^2\beta^2 \mathcal{T}x^2y^2+4\alpha\beta ^3\mathcal{T}x y^3+\beta^4 \mathcal{T}y^4,\endaligned$$
$$\mathcal{A}(\alpha x+\beta y)^3=\langle\mathcal{A}, (\alpha x+\beta y)^{\otimes 3}\rangle=\alpha^4\mathcal{A}x^3+3\alpha^2\beta \mathcal{A}x^2 y
+3\alpha\beta ^2\mathcal{A}x y^2+\beta^3 \mathcal{A}y^3$$
where $\mathcal{T}x^ky^l=\sum\limits_{i_1,i_2, \cdots, i_d = 1}^nt_{i_1i_2\cdots i_d}x_{i_1}x_{i_2}\cdots x_{i_k} y_{i_{k+1}}y_{i_{k+2}}\cdots y_{i_{k+l}}$ and $k+l=d.$

\section{Positive definiteness of 4th order 3-dimensional symmetric tensors}
A 4th order 3-dimensional  tensor has 81 components, 
$$\begin{pmatrix}
		\begin{pmatrix}
			t_{1111} & t_{1211} & t_{1311}\\
			t_{2111} & t_{2211} & t_{2311}\\
			t_{3111} & t_{3211} & t_{3311}
		\end{pmatrix}& \begin{pmatrix}
			t_{1112} & t_{1212} & t_{1312}\\
			t_{2112} & t_{2212} & t_{2312}\\
			t_{3112} & t_{3212} & t_{3312}
		\end{pmatrix}& \begin{pmatrix}
			t_{1113} & t_{1213} & t_{1313}\\
			t_{2113} & t_{2213} & t_{2313}\\
			t_{3113} & t_{3213} & t_{3313}
		\end{pmatrix}\\
		& & \\
		\begin{pmatrix}
			t_{1121} & t_{1221} & t_{1321}\\
			t_{2121} & t_{2221} & t_{2321}\\
			t_{3121} & t_{3221} & t_{3321}
		\end{pmatrix}& \begin{pmatrix}
			t_{1122} & t_{1222} & t_{1322}\\
			t_{2122} & t_{2222} & t_{2322}\\
			t_{3122} & t_{3222} & t_{3322}
		\end{pmatrix}& \begin{pmatrix}
			t_{1123} & t_{1223} & t_{1323}\\
			t_{2123} & t_{2223} & t_{2323}\\
			t_{3123} & t_{3223} & t_{3323}
		\end{pmatrix}\\
		& & \\
		\begin{pmatrix}
			t_{1131} & t_{1231} & t_{1331}\\
			t_{2131} & t_{2231} & t_{2331}\\
			t_{3131} & t_{3231} & t_{3331}
		\end{pmatrix}& \begin{pmatrix}
			t_{1132} & t_{1232} & t_{1332}\\
			t_{2132} & t_{2232} & t_{2332}\\
			t_{3132} & t_{3232} & t_{3332}
		\end{pmatrix}& \begin{pmatrix}
			t_{1133} & t_{1233} & t_{1333}\\
			t_{2133} & t_{2233} & t_{2333}\\
			t_{3133} & t_{3233} & t_{3333}
		\end{pmatrix}
\end{pmatrix}$$
Moreover, such a tensor has only 15 free entries if it is symmetric, 
$$\aligned &t_{1111},  t_{2222}, t_{3333}, t_{1122}, t_{1133}, t_{2233}, t_{1112}, t_{1113}, \\&t_{1222}, t_{1333}, t_{2223}, t_{2333}, t_{1123}, t_{1223}, t_{1233}.\endaligned$$
Then for $x=(x_1,x_2,x_3)^\top$, a ternary quartic form (homogeneous polynomial) is uniquely given by a fourth order symmetric tensor $\mathcal{T}$, 
$$\aligned 
\mathcal{T}x^4=&t_{1111}x_1^4+t_{2222}x_2^4+t_{3333}x_3^4+6t_{1122}x_1^2x_2^2+6t_{1133}x_1^2x_3^2+6t_{2233}x_2^2x_3^2\\&+4t_{1112}x_1^3x_2+4t_{1113}x_1^3x_3+4t_{1222}x_1x_2^3+4t_{1333}x_1x_3^3+4t_{2223}x_2^3x_3+4t_{2333}x_2x_3^3\\
&+12t_{1123}x_1^2x_2x_3+12t_{1223}x_1x_2^2x_3+12t_{1233}x_1x_2x_3^2.
\endaligned$$
\begin{definition}\label{def:14} A fourth order $3$-dimensional   tensor $\mathcal{T}=(t_{ijkl})$ is said to be
		\begin{itemize}
			\item {\bf cyclic symmetric} if it is symmetric and $$\aligned &t_{1111}=t_{2222}=t_{3333},		t_{1112}=t_{2223}=t_{1333}, t_{1113}=t_{1222}=t_{2333},\\
			& t_{1122}=t_{1133}=t_{2233},  t_{1123}=t_{1223}=t_{1233};\endaligned $$
			\item {\bf totally symmetric} if it is symmetric and $$\aligned &t_{1111}=t_{2222}=t_{3333}, t_{1112}=t_{2223}=t_{1333}= t_{1113}=t_{1222}=t_{2333},\\
			& t_{1122}=t_{1133}=t_{2233},  t_{1123}=t_{1223}=t_{1233};\endaligned $$
		\end{itemize} 
		
		\end{definition}
\begin{theorem}
	\label{thm:41} \em Let $\mathcal{T}=(t_{ijkl})$ be a 4th-order $n$-dimensional symmetric tensor. Then 
	$\mathcal{T}$ is positive definite if and only if $$\begin{cases}
		\mathcal{T}x^4=0\ \Longrightarrow\ x=0,\\
		\mbox{there is a }y\in\mathbb{R}^n\setminus\{0\}\mbox{ such that }\mathcal{T}y^4>0;
	\end{cases}$$
\end{theorem}
{\bf Proof.}
	The necessarity is obvious.  Now we show the sufficiency. Suppose $\mathcal{T}$ is not positive definite when the conditions are satisfied. There exists $u\in\mathbb{R}^n\setminus \{0\}$ such that $\mathcal{T}u^4\leq0$.  Since $\mathcal{T}u^4=0$ means $u=0$ by the conditions, then $\mathcal{T}u^4<0$. Apply the intermediate value theore to continuous function $\mathcal{T}x^4$, there is an $\lambda\in(0,1)$ such that \begin{center}
		$z=(1-\lambda)u+\lambda y$ satisfying $\mathcal{T}z^4 =0.$
	\end{center} This implies $z=0$, i.e., $y=\dfrac{\lambda-1}{\lambda}u$. So, we find
	$$\aligned\mathcal{T}z^4=&\left\langle\mathcal{T}, ((1-\lambda)u+\lambda y)^{\otimes 4}\right\rangle=(1-\lambda)^4\mathcal{T}u^4+4(1-\lambda)^3\lambda \mathcal{T}u^3 y\\
	&+6(1-\lambda)^2\lambda^2 \mathcal{T}u^2y^2+4(1-\lambda)\lambda ^3\mathcal{T}u y^3+\lambda^4 \mathcal{T}y^4\\
	=&(1-\lambda)^4\mathcal{T}u^4-4(1-\lambda)^4 \mathcal{T}u^4 
	+6(1-\lambda)^4 \mathcal{T}u^4-4(1-\lambda)^4\mathcal{T}u^4+\lambda^4 \mathcal{T}y^4\\
	=&-(1-\lambda)^4 \mathcal{T}u^4+\lambda^4 \mathcal{T}y^4>0,\endaligned$$
	a contradiction. Therefore, $\mathcal{T}$ is  positive definite. 
	
\begin{theorem}\label{thm:42} \em Let $\mathcal{T}=(t_{ijkl})$ be a 4th-order $3$-dimensional cyclic symmetric tensor. Suppose $\mathcal{T}$ is  positive semi-definite and \begin{center}
		$|t_{1112}|=|t_{1222}|=t_{1122}=t_{1111}=1$ and $t_{1112}t_{1222}=-1$.
	\end{center}
	Then  $t_{1123}\ge-\dfrac7{12}.$
\end{theorem}
{\bf Proof.} By the cyclic symmetry of $\mathcal{T}$, we might take $t_{1222}=t_{2333}=t_{1113}=1$ and $t_{1112}=t_{1333}=t_{2223}=-1$ without loss the generality. 

  From the  positive semidefiniteness of $\mathcal{T}$, it follows that for $x=(1,1,1)^\top$, 
$$\aligned 
\mathcal{T}x^4=&x_1^4+x_2^4+x_3^4+6(x_1^2x_2^2+x_1^2x_3^2+x_2^2x_3^2)\\&+4(x_1x_2^3+x_1^3x_3+x_2x_3^3)-4(x_1^3x_2+x_1x_3^3+x_2^3x_3)\\
&+12t_{1123}(x_1^2x_2x_3+x_1x_2^2x_3+x_1x_2x_3^2)\\
=&(x_1+x_2+x_3)^4-8(x_1^3x_2+x_1x_3^3+x_2^3x_3)\\&+12(t_{1123}-1)x_1x_2x_3(x_1+x_2+x_3)\\
=&3^4-8\times 3+36(t_{1123}-1)\geq0,\endaligned$$
and hence, $$
t_{1123}\ge -\dfrac7{12}.
$$

\begin{theorem}\label{thm:43} \em Let $\mathcal{T}=(t_{ijkl})$ be a 4th-order $3$-dimensional cyclic symmetric tensor. Suppose  \begin{center}
		$|t_{1112}|=|t_{1222}|=t_{1122}=t_{1111}=1$ and $t_{1123}=-\dfrac7{12}$.  
	\end{center}
	Then  $\mathcal{T}$ is  positive semi-definite if and only if $t_{1112}t_{1222}=-1$.
\end{theorem}
{\bf Proof.}  ``only if. " Suupose $\mathcal{T}$ is  positive semi-definite, but the condition can't hold. Then $t_{1112}t_{1222}=1$. That is, $t_{1112}=t_{1222}=1$ or $t_{1112}=t_{1222}=-1$. 

Case 1. $t_{1222}=t_{2333}=t_{1113}=t_{1112}=t_{1333}=t_{2223}=1$.   Then for  $x=(1,1,-5)^\top$, 
$$\aligned 
\mathcal{T}x^4=&(x_1+x_2+x_3)^4+12(t_{1123}-1)x_1x_2x_3(x_1+x_2+x_3)\\
=&(x_1+x_2+x_3)^4-19x_1x_2x_3(x_1+x_2+x_3)\\
=&(1+1-5)^4+19\times5(1+1-5)=-204<0.
\endaligned$$ This is a contradiction to  the positive semidefiniteness of $\mathcal{T}$.

Case 2. $t_{1222}=t_{2333}=t_{1113}=t_{1112}=t_{1333}=t_{2223}=-1$. For $x=(1,1,1)^\top$, 
$$\aligned 
\mathcal{T}x^4=&(x_1+x_2+x_3)^4-8(x_1^3x_2+x_1x_3^3+x_2^3x_3+x_1x_2^3+x_1^3x_3+x_2x_3^3)\\&-19x_1x_2x_3(x_1+x_2+x_3)\\
=&3^4-8\times 6-19\times 3=-24<0,\endaligned$$
a contradiction.  So with all that, the condition, $t_{1112}t_{1222}=-1$ is necessary.

``if. " Assume $t_{1112}t_{1222}=-1$.  Then we might take $t_{1222}=t_{2333}=t_{1113}=1$ and $t_{1112}=t_{1333}=t_{2223}=-1$. Rewriting $\mathcal{T}x^4$, $$\aligned \mathcal{T}x^4=&(x_1+x_2+x_3)^4-8(x_1^3x_2+x_1x_3^3+x_2^3x_3)-19x_1x_2x_3(x_1+x_2+x_3)\\
=&(x_1-x_2)^4+(x_1-x_3)^4+(x_3-x_2)^4-x_1^4-x_2^4-x_3^4\\
&+8(x_1x_2^3+x_1^3x_3+x_2x_3^3)-7x_1x_2x_3(x_1+x_2+x_3).\endaligned$$
Solve the optimization problem:
$$\min\{\mathcal{T}x^4: x_1^2+x_2^2+x_3^2=1\}$$
to yield  the minimum value $0$ of $\mathcal{T}x^4$ at two minimizers $$\left(\dfrac1{\sqrt3},\dfrac1{\sqrt3},\dfrac1{\sqrt3}\right)^\top,\ \  \left(-\dfrac1{\sqrt3},-\dfrac1{\sqrt3},-\dfrac1{\sqrt3}\right)^\top.$$
So, $\mathcal{T}$ is  positive semi-definite.

\begin{theorem}\label{thm:44} \em Let $\mathcal{T}=(t_{ijkl})$ be a 4th-order $3$-dimensional cyclic symmetric tensor. Suppose \begin{center}
	$t_{1123}\in\big(-\frac7{12}, -\frac16 \big]$,	$|t_{1112}|=|t_{1222}|=t_{1122}=t_{1111}=1$ and $t_{1112}t_{1222}=-1$.
	\end{center}
	 Then  $\mathcal{T}$ is  positive definite.
\end{theorem}
{\bf Proof.} Let $t_{1222}=t_{2333}=t_{1113}=1$ and $t_{1112}=t_{1333}=t_{2223}=-1$ without loss the generality.  For  $t_{1123}=t_{1223}=t_{1233}\in\big(-\frac7{12}, -\frac16 \big]$,  $\mathcal{T}x^4$ may be rewritten as
	\begin{equation}\label{eq:41}\aligned
		\mathcal{T}x^4=&(x_1+x_2+x_3)^4-8(x_1^3x_2+x_1x_3^3+x_2^3x_3)\\&+12(t_{1123}-1)x_1x_2x_3(x_1+x_2+x_3).\endaligned
	\end{equation} 
	Then it can be divided into three kinds of circumstances.
	
	\textbf{Case 1.} $x_1x_2x_3(x_1+x_2+x_3)=0$.  This implies that $x_1=0$ or $x_2=0$ or $x_3=0$ or $x_1+x_2+x_3=0$, and then, the top three conditions definited respectively a 4th order 2-dimensional principal subtensor. Therefore, an immediate application of Corollary \ref{cor:23} yields $\mathcal{T}x^4>0$ for all $x\ne0 $ under the top three conditions. Next we show that $\mathcal{T}x^4>0$ for all $x\ne0 $ with 
	$x_1=-x_2-x_3$. In fact, by this time,  we have $$\aligned \mathcal{T}x^4=&-8(-(x_2+x_3)^3x_2-(x_2+x_3)x_3^3+x_2^3x_3)\\
	=&8(x_2^4+2x_2^3x_3+3x^2_2x_3^2+2x_2x_3^3+x_3^4)\\
	=&8\mathcal{A}x^4,
	\endaligned$$
	where $\mathcal{A}=(a_{ijkl})$  is a 4th order 2-dimensional tensor with its entries
	$$a_{2222}=a_{3333}=1, a_{2223}=a_{2333}= a_{2233}=\frac12,$$
	and then,  such a  tensor $\mathcal{A}$ satisfies the conditions of Corollary \ref{cor:21}, $$a_{2223}=\dfrac12, 2a_{2223}^2+1=3a_{2233}, a_{2223}=a_{2333}.$$ 
	So an application of Corollary \ref{cor:21} yields that $\mathcal{T}x^4>0$ for all $x\ne0 $ with $x_1=-x_2-x_3$.
	
	\textbf{Case 2.} $x_1x_2x_3(x_1+x_2+x_3)>0$.  That's when $\mathcal{T}x^4$ is always greater than the value of $\mathcal{T}x^4$ with $t_{1123}=-\frac7{12}$. That is, $$\aligned \mathcal{T}x^4=&(x_1+x_2+x_3)^4-8(x_1^3x_2+x_1x_3^3+x_2^3x_3)+12(t_{1123}-1)x_1x_2x_3(x_1+x_2+x_3)
	\\
	>&(x_1+x_2+x_3)^4-8(x_1^3x_2+x_1x_3^3+x_2^3x_3)-19x_1x_2x_3(x_1+x_2+x_3), \endaligned$$
	and so, it follows from Theorem \ref{thm:42} that $\mathcal{T}x^4>0$ under this condition.
	
	\textbf{Case 3. } $x_1x_2x_3(x_1+x_2+x_3)<0$.  By this time, for each $x$,  $\mathcal{T}x^4$ always reaches its minimum value at $t_{1123}=-\frac16$. That is, $$ \mathcal{T}x^4
	\geq(x_1+x_2+x_3)^4-8(x_1^3x_2+x_1x_3^3+x_2^3x_3)-14x_1x_2x_3(x_1+x_2+x_3).$$
	Let  $\mathcal{T}'x^4$ be definited by $$\aligned \mathcal{T}'x^4=&(x_1+x_2+x_3)^4-8(x_1^3x_2+x_1x_3^3+x_2^3x_3)-14x_1x_2x_3(x_1+x_2+x_3).\endaligned$$
	We  directly solve the equation, $\mathcal{T}'x^4=0$ to yield $x=0$, and then apply Theorem \ref{thm:41} to establish the conclusion that $\mathcal{T}'$ is positive definite. Therefore, $\mathcal{T}x^4>0$ under this condition.
	
	So with all that, $\mathcal{T}$ is positive definite.

\begin{corollary}\label{cor:45}
	\em Let $\mathcal{T}=(t_{ijkl})$ be a 4th-order $3$-dimensional cyclic symmetric tensor.  Assume that $$t_{1111}=|t_{1112}|=|t_{1222}|=1, t_{1112}t_{1222}=-1.$$ Then 
	\begin{itemize}
		\item [(i)]\ \  $\mathcal{T}$ is positive semi-definite if  $t_{1123}\in\big[-\frac7{12}, -\frac16 \big]$ and $t_{1122}\ge1$;
		\item [(ii)]\ \  $\mathcal{T}$ is positive definite   if $t_{1123}\in\big(-\frac7{12}, -\frac16 \big]$  and $t_{1122}\ge1$.
	\end{itemize}
\end{corollary}
{\bf Proof.}  Since $t_{1122}\ge1$ and $\mathcal{T}$ is cyclic symmetric,  for each $x=(x_1,x_2,x_3)^\top\in\mathbb{R}^3$, we have $$\aligned 
	\mathcal{T}x^4=&x_1^4+x_2^4+x_3^4+4t_{1112}(x_1^3x_2+x_1x_3^3+x_3x_2^3)+4t_{1222}(x_1^3x_3+x_1x_2^3+x_2x_3^3)\\&+6t_{1122}(x_1^2x_2^2+x_1^2x_3^2+x_2^2x_3^2)+12t_{1123}x_1x_2x_3(x_1+x_2+x_3)\\
	\ge&x_1^4+x_2^4+x_3^4+4t_{1112}(x_1^3x_2+x_1x_3^3+x_3x_2^3)+4t_{1222}(x_1^3x_3+x_1x_2^3+x_2x_3^3)\\&+6(x_1^2x_2^2+x_1^2x_3^2+x_2^2x_3^2)+12t_{1123}x_1x_2x_3(x_1+x_2+x_3)\\
	=&	\mathcal{T}'x^4.
	\endaligned$$ Then a tensor  $\mathcal{T}'$ reaches the conditions of  Theorems \ref{thm:43} ($t_{1123}=-\frac7{12}$)  and \ref{thm:44}, and so,  the desired results follow.\\
 
 The positive semidefiniteness of $\mathcal{T}$ does not have to require  three elements, $t_{1123}= t_{1223}=t_{1233}$.

\begin{theorem}\label{thm:46} \em Let $\mathcal{T}=(t_{ijkl})$ be a 4th-order $3$-dimensional symmetric tensor.   Assume $t_{1222}=t_{2333}=t_{1113}$, $t_{1112}=t_{1333}=t_{2223}$ and \begin{center}
		$|t_{iiij}|=t_{iiii}=t_{iijj}=1$ and $t_{ijjj}t_{iiij}=-1$ for all $i,j\in\{1,2,3\}$, $i\ne j$.
	\end{center}
	Then  \begin{itemize}
		\item[(i)] \ \  $\mathcal{T}$ is  positive definite if	$t_{1123}, t_{1223}, t_{1233}\in\big(-\frac7{12}, -\frac14 \big]$;
		\item[(ii)] \ \  $\mathcal{T}$ is  positive definite if	 $t_{1123}, t_{1223}, t_{1233}\in\big[-\frac14,-\frac16 \big]$.
	\end{itemize}
\end{theorem}
{\bf Proof.}
	Let $t_{1222}=t_{2333}=t_{1113}=1$ and $t_{1112}=t_{1333}=t_{2223}=-1$ without loss the generality. 
	
	(i) For all $t_{1123}, t_{1223}, t_{1233}\in\big(-\frac7{12}, -\frac14 \big]$,  $\mathcal{T}x^4$ may be rewritten as
	\begin{equation}\label{eq:42}\aligned
		\mathcal{T}x^4=&(x_1+x_2+x_3)^4-8(x_1^3x_2+x_1x_3^3+x_2^3x_3)\\&+12(t_{1123}-1)x_1^2x_2x_3+12(t_{1223}-1)x_1x_2^2x_3\\&+12(t_{1233}-1)x_1x_2x_3^2.
		\endaligned
\end{equation} 
	Then it can be divided into two cases.
	
	\textbf{Case 1.} $x_1^2x_2x_3=0$ or $x_1x_2^2x_3=0$ or $x_1x_2x_3^2=0$.  This means that $x_1=0$ or $x_2=0$ or $x_3=0$, which respectively gives a 4th order 2-dimensional principal subtensor, and hence, an immediate application of Corollary \ref{cor:23} yields that $\mathcal{T}x^4>0$ for all $x\ne0 $ under this case. 
	
	\textbf{Case 2.} $x_1^2x_2x_3\ne0$ and $x_1x_2^2x_3\ne0$ and $x_1x_2x_3^2\ne0$, which implies $x_1\ne0$ and $x_2\ne0$ and $x_3\ne0$.  This can be divided into three subcases.
	
	\textbf{Subcase 1.} $t_{1123}=t_{1223}=t_{1233}$.  The conclusion directly follows from Theorem \ref{thm:44}.  
	
	\textbf{Subcase 2.} Two of three elements $t_{1123}, t_{1223}, t_{1233}$ are equal. We might take $t_{1223}=t_{1233}$ without loss of generality.    Rewrite Eq. \eqref{eq:42}  as
	\begin{equation}\label{eq:43}\aligned
		\mathcal{T}x^4=&(x_1+x_2+x_3)^4-8(x_1^3x_2+x_1x_3^3+x_2^3x_3)\\&+12(t_{1123}-1)x_1^2x_2x_3+12(t_{1223}-1)x_1x_2x_3(x_2+x_3).\endaligned
	\end{equation} 
	If $x_2+x_3=0$, then we have 
	$$\aligned
	\mathcal{T}x^4=&x_1^4-8(x_1^3x_2-x_1x_2^3-x_2^3x_2)-12(t_{1123}-1)x_1^2x_2^2\\
	=&x_1^4-8x_1^3x_2+8x_1x_2^3+8x_2^4-12(t_{1123}-1)x_1^2x_2^2\\
	\geq& x_1^4-8x_1^3x_2+15x_1^2x_2^2+8x_1x_2^3+8x_2^4.\endaligned$$
	Let $\mathcal{B}x^4=x_1^4-8x_1^3x_2+15x_1^2x_2^2+8x_1x_2^3+8x_2^4.$ This gives a 4th order 2-dimensional  tensor $\mathcal{B}=(b_{ijkl})$ with its entries
	$$b_{1111}=1, b_{2222}=8, b_{1112}=-2, b_{1222}=2, b_{1122}=\frac52.$$
	Then by Lemma \ref{lem:21}, $\mathcal{B}$ is positive definite, and so, $\mathcal{T}x^4>0$ for all $x\in\mathbb{R}^3\setminus\{0\}$ with $x_2+x_3=0.$
	
	If $x_2+x_3\ne0$, then  we have $$\begin{cases}
		12(t_{1123}-1)x_1^2x_2x_3>-19x_1^2x_2x_3, &  x_2x_3>0;\\
		12(t_{1123}-1)x_1^2x_2x_3\geq-15x_1^2x_2x_3, &  x_2x_3<0,
	\end{cases}$$
	and 
	$$\begin{cases}
		12(t_{1223}-1)x_1x_2x_3(x_2+x_3)>-19x_1x_2x_3(x_2+x_3), &  x_1x_2x_3(x_2+x_3)>0;\\
		12(t_{1223}-1)x_1x_2x_3(x_2+x_3)\geq-15x_1x_2x_3(x_2+x_3), &  x_1x_2x_3(x_2+x_3)<0.
	\end{cases}$$
	Thus, from Eq. \eqref{eq:43}, they follow  that when $x_1^2x_2x_3<0$ and $x_1x_2x_3(x_2+x_3)<0$, 
	$$\mathcal{T}x^4\geq (x_1+x_2+x_3)^4-8(x_1^3x_2+x_1x_3^3+x_2^3x_3)-x_1x_2x_3(15x_1+15x_2+15x_3);$$ when  $x_1^2x_2x_3>0$ and $x_1x_2x_3(x_2+x_3)<0$, 
	$$\mathcal{T}x^4> (x_1+x_2+x_3)^4-8(x_1^3x_2+x_1x_3^3+x_2^3x_3)-x_1x_2x_3(19x_1+15x_2+15x_3);$$ 
	when $x_1^2x_2x_3<0$ and $x_1x_2x_3(x_2+x_3)>0$,
	$$\mathcal{T}x^4> (x_1+x_2+x_3)^4-8(x_1^3x_2+x_1x_3^3+x_2^3x_3)-x_1x_2x_3(15x_1+19x_2+19x_3);$$
	when $x_1^2x_2x_3>0$ and $x_1x_2x_3(x_2+x_3)>0$,  
	$$\mathcal{T}x^4> (x_1+x_2+x_3)^4-8(x_1^3x_2+x_1x_3^3+x_2^3x_3)-x_1x_2x_3(19x_1+19x_2+19x_3).$$
	Similarly to Case 3 of Theorem \ref{thm:44}, 
	it is easy to verify that in  the above right polynomials, the top three ones  are all positive and the fourth one is non-negative  for all $x\ne0$.  Therefore, $\mathcal{T}x^4>0$ for all $x\in\mathbb{R}^3\setminus\{0\}$  with $x_2+x_3\ne0$.
	
	\textbf{Subcase 3.}   $t_{1123}\ne t_{1223}\ne t_{1233}$, which are mutually different.  For all $t_{1123}, t_{1223}, t_{1233}\in\big(-\frac7{12}, -\frac14 \big]$, 
	similarly to \textbf{Subcase 2}, we also have 
	$$\begin{cases}
		12(t_{1223}-1)x_1x_2^2x_3>-19x_1x_2^2x_3, & x_1x_3>0;\\
		12(t_{1223}-1)x_1x_2^2x_3\geq-15x_1x_2^2x_3, & x_1x_3<0,
	\end{cases}$$ and 
	$$\begin{cases}
		12(t_{1233}-1)x_1x_2x_3^2>-19x_1x_2x_3^2, & x_1x_2>0;\\
		12(t_{1233}-1)x_1x_2x_3^2\geq-15x_1x_2x_3^2, & x_1x_2<0.
	\end{cases}$$
	However,  three inequalities, $x_2x_3<0$, $x_1x_3<0$ and $x_1x_2<0$ can't simultaneously hold. In fact,  the inequalities $x_1x_3<0$ and $x_1x_2<0$ imply that $x_3$ and $x_2$ have the same sign, and hence, $x_2x_3>0$. 
	
	Three inequalities, $x_2x_3<0$, $x_1x_3>0$ and $x_1x_2>0$ are incompatible. In fact, two inequalities $x_2x_3<0$ and $x_1x_3>0$ mean that $x_2$ and $x_1$ 
	have opposite sign, and hence, $x_1x_2<0$.  So there doesn't exist this case: two of three numbers, $x_2x_3$, $x_1x_3$ and $x_1x_2$ are positive, and one of them is negative.
	
	Thus, from Eq. \eqref{eq:42}, it follows that  when two of $x_1^2x_2x_3$, $x_1x_2^2x_3$ and $x_1x_2x_3^2$ is negative,  without loss of generality, we may take $x_1^2x_2x_3>0$ and $x_1x_2^2x_3<0$ and $x_1x_2x_3^2<0$, and then
	$$\mathcal{T}x^4> (x_1+x_2+x_3)^4-8(x_1^3x_2+x_1x_3^3+x_2^3x_3)-x_1x_2x_3(19x_1+15x_2+15x_3);$$ 
	when $x_1^2x_2x_3>0$ and $x_1x_2^2x_3>0$ and $x_1x_2x_3^2>0$,  
	$$\mathcal{T}x^4> (x_1+x_2+x_3)^4-8(x_1^3x_2+x_1x_3^3+x_2^3x_3)-x_1x_2x_3(19x_1+19x_2+19x_3).$$
	Similarly to Case 3 of Theorem \ref{thm:44}, 
	it is easy to verify that in   the above right  polynomials,  the first one  is  positive and the second one is non-negative for all $x\ne0$.  Therefore, $\mathcal{T}x^4>0$ for all $x\ne0$ under this subcase.
	
	So with all that, $\mathcal{T}$ is positive definite.
	
	(ii) Using the same argumentation technique of (i), the desired conclusion follows. This completes the proof.

\begin{remark}\label{rem:41}\em
	In Theorem \ref{thm:46}, two intervals $\big(-\frac7{12}, -\frac14 \big]$ and $\big[-\frac14,-\frac16 \big]$ can't write as their union, $\big(-\frac7{12}, -\frac16 \big]$ because  the four polynomials 
	$$\aligned (x_1+x_2+x_3)^4-8(x_1^3x_2+x_1x_3^3+x_2^3x_3)-x_1x_2x_3(19x_1+14x_2+14x_3),\\
	(x_1+x_2+x_3)^4-8(x_1^3x_2+x_1x_3^3+x_2^3x_3)-x_1x_2x_3(18x_1+14x_2+14x_3),\\
	(x_1+x_2+x_3)^4-8(x_1^3x_2+x_1x_3^3+x_2^3x_3)-x_1x_2x_3(17x_1+14x_2+14x_3),\\
	(x_1+x_2+x_3)^4-8(x_1^3x_2+x_1x_3^3+x_2^3x_3)-x_1x_2x_3(16x_1+14x_2+14x_3),\\
	\endaligned$$
	are not all positive for all $x\in\mathbb{R}^3\setminus\{0\}$ (for example, $x=(-1.2,5,1)^\top$). But the polynimials, $$(x_1+x_2+x_3)^4-8(x_1^3x_2+x_1x_3^3+x_2^3x_3)-x_1x_2x_3(15x_1+14x_2+14x_3)>0$$ for all $x\in\mathbb{R}^3\setminus\{0\}$ and $$(x_1+x_2+x_3)^4-8(x_1^3x_2+x_1x_3^3+x_2^3x_3)-x_1x_2x_3\left(\dfrac{46}3x_1+14x_2+14x_3\right)>0$$ for all $x\in\mathbb{R}^3\setminus\{0\}$.  Therefore, if we respectively replace two intervals $\big(-\frac7{12}, -\frac14 \big]$ and $\big[-\frac14,-\frac16 \big]$ with $\big(-\frac7{12}, -\frac5{18} \big]$ and $\big[-\frac5{18},-\frac16 \big]$  in Theorem \ref{thm:46},  then the conclusions still hold.\end{remark}
	
	\begin{remark}\label{rem:42}\em
	In Theorem \ref{thm:44}, the right node $-\frac16 $ of the interval $\big(-\frac7{12} , -\frac16 \big]$ may be replaced by $-\frac5{36}$ since it is easy to verify a polynomial $$(x_1+x_2+x_3)^4-8(x_1^3x_2+x_1x_3^3+x_2^3x_3)-\dfrac{41}3x_1x_2x_3(x_1+x_2+x_3)>0$$ for all $x\in\mathbb{R}^3\setminus\{0\}.$ So, in Theorem \ref{thm:44} and  Corollary \ref{cor:45}, the interval $\big(-\frac7{12} , -\frac5{36} \big]$ can substitute for $\big(-\frac7{12} , -\frac16 \big]$. However, in Theorem \ref{thm:46}, $\big[-\frac14,-\frac5{36} \big]$ can not substitute for $\big[-\frac14,-\frac16 \big]$ as a polynomial $$(x_1+x_2+x_3)^4-8(x_1^3x_2+x_1x_3^3+x_2^3x_3)-x_1x_2x_3\left(\dfrac{41}3x_1+15x_2+15x_3\right)>0$$ may not hold  for all $x\in\mathbb{R}^3\setminus\{0\}$ (for example, $x=(-9.4,-2,2.3)^\top)$.
\end{remark}

By applying Theorems \ref{thm:42}, \ref{thm:44}  and \ref{thm:46},  the following  inequalities are established easily  for  ternary quartic homogeneous polynomials.

\begin{corollary}\label{cor:47}	If $(x_1,x_2,x_3)\ne(0,0,0)$, then 
	\begin{itemize}
		\item [(i)] $(x_1+x_2+x_3)^4\geq8(x_1^3x_2+x_1x_3^3+x_2^3x_3)+19x_1x_2x_3(x_1+x_2+x_3),
		$\\ with equality if and only if  $x_1= x_2=x_3$;
		\item [(ii)] $(x_1+x_2+x_3)^4>8(x_1^3x_2+x_1x_3^3+x_2^3x_3)+14x_1x_2x_3(x_1+x_2+x_3);$
		\item [(iii)] $(x_1+x_2+x_3)^4>8(x_1^3x_2+x_1x_3^3+x_2^3x_3)+15x_1x_2x_3(x_1+x_2+x_3);$
		\item [(iv)] $(x_1+x_2+x_3)^4>8(x_1^3x_2+x_1x_3^3+x_2^3x_3)+16x_1x_2x_3(x_1+x_2+x_3);$
		\item [(v)] $(x_1+x_2+x_3)^4>8(x_1^3x_2+x_1x_3^3+x_2^3x_3)+17x_1x_2x_3(x_1+x_2+x_3);$
		\item [(vi)] $(x_1+x_2+x_3)^4>8(x_1^3x_2+x_1x_3^3+x_2^3x_3)+18x_1x_2x_3(x_1+x_2+x_3).$
	\end{itemize} 
	Furthermore, these (strict) inequalities still hold if $x_1^3x_2$ and $x_1x_2^3$,  $x_1^3x_3$ and $x_1x_3^3$, $x_2x_3^3$ and  $x_2^3x_3$ are simultaneously  exchangeable.
\end{corollary}
	\begin{corollary}\label{cor:48}	If $(x_1,x_2,x_3)\ne(0,0,0)$, then 
	\begin{itemize}
		\item [(i)] $(x_1+x_2+x_3)^4>8(x_1^3x_2+x_1x_3^3+x_2^3x_3)+x_1x_2x_3(19x_1+17x_2+15x_3);$
		\item [(iii)] $(x_1+x_2+x_3)^4>8(x_1^3x_2+x_1x_3^3+x_2^3x_3)+x_1x_2x_3(19x_1+16x_2+15x_3);$
		\item [(iv)] $(x_1+x_2+x_3)^4>8(x_1^3x_2+x_1x_3^3+x_2^3x_3)+x_1x_2x_3(15x_1+14x_2+14x_3);$
		\item [(v)] $(x_1+x_2+x_3)^4>8(x_1^3x_2+x_1x_3^3+x_2^3x_3)+x_1x_2x_3(15x_1+16x_2+14x_3);$
		\item [(vi)] $(x_1+x_2+x_3)^4>8(x_1^3x_2+x_1x_3^3+x_2^3x_3)+x_1x_2x_3(17x_1+15x_2+18x_3).$
	\end{itemize} 
	Furthermore, these strict inequalities still hold if $x_1^3x_2$ and $x_1x_2^3$,  $x_1^3x_3$ and $x_1x_3^3$, $x_2x_3^3$ and  $x_2^3x_3$ are simultaneously  exchangeable.
\end{corollary}

\section{Conclusions}

For a class of 4th order 3-dimensional  tensors,  the conditions testing its positive definiteness  are obtained under such  tensors with cyclic symmetry and  non-cyclic symmetry, respectively.  The necessary and sufficient conditions are established for the positive (semi-)definiteness of fourth order tensors.    Several  (strict) inequalities  of quartic ternary  homogeneous polynomial are showed with the help of these analytic conditions.

\end{document}